Testing Surrogate-Based Optimization with the Fortified Branin-Hoo Extended to Four Dimensions

Charles F. Jekel[a,b,1] and Raphael T. Haftka[a]

a. Work performed at Department of Mechanical and Aerospace Engineering, University of Florida.
b. Currently at Lawrence Livermore National Laboratory[2], Livermore, CA.

**Abstract:** Some popular functions used to test global optimization algorithms have multiple local optima, all with the same value, making them all global optima. It is easy to make them more challenging by fortifying them via adding a localized bump at the location of one of the optima. In previous work the authors illustrated this for the Branin-Hoo function and the popular differential evolution algorithm, showing that the fortified Branin-Hoo required an order of magnitude more function evaluations. This paper examines the effect of fortifying the Branin-Hoo function on surrogate-based optimization, which usually proceeds by adaptive sampling. Two algorithms are considered. The EGO algorithm, which is based on a Gaussian process (GP) and an algorithm based on radial basis functions (RBF). EGO is found to be more frugal in terms of the number of required function evaluations required to identify the correct basin, but it is expensive to run on a desktop, limiting the number of times the runs could be repeated to establish sound statistics on the number of required function evaluations. The RBF algorithm was cheaper to run, providing more sound statistics on performance. A four-dimensional version of the Branin-Hoo function was introduced in order to assess the effect of dimensionality. It was found that the difference between the ordinary function and the fortified one was much more pronounced for the four-dimensional function compared to the two dimensional one.

I.   **Introduction:**

Some popular functions used to test global optimization algorithms have multiple local optima, all with the same value, making them all global optima. The Branin-Hoo function is a typical example, having three global optima. Jekel and Haftka (2019, 2020) demonstrated that it is easy to fortify such functions by adding a localized bump at the location of one of the optima, making them much more challenging to optimize in terms of the required number of function evaluations. This was illustrated for the Branin-Hoo function for the popular differential evolution (DE, Storn and Price, 1997) algorithm. It was found that the fortified Branin-Hoo function required an order of magnitude more function evaluations.

Surrogate-based optimization is considered to be more frugal in terms of the required number of function evaluations than algorithms such as DE that do not use a surrogate. Surrogate-based optimization strategies typically take advantage of adaptive sampling (Forrester and Keane, 2009). That is, a surrogate is fitted to an initial sample of points, and then some algebraic objective is optimized to locate the next sampling point. The surrogate is refitted with the extra point, and the process is repeated until a stopping criterion is met.

The efficient global optimization (EGO, Jones et al. 1998) algorithm uses a Gaussian process (GP) or kriging surrogate for that purpose, selecting the next sampling point by minimizing the expected improvement (EI) calculated based on the uncertainty model provided by the GP. The algorithm is substantially more expensive to run than DE, as fitting a GP and maximizing EI is more expensive than the crossover and mutation calculation of DE. On a desktop computer, the runtime for an EGO iteration could be seconds compared to microseconds of a DE iteration (when considering a trivial analytical objective function). However, challenging engineering simulations often require days or weeks on expensive computers using thousands of processors to compute an objective function, so the additional computational overhead of EGO itself is negligible. However, thousands of replicate runs are needed in order to compare the performance of the algorithm on the original and fortified functions. This is because both DE and EGO have randomness built into them so that the final results also have some random element.

Because EGO did not permit us testing with thousands of repeated runs, we also considered an algorithm based on regression with radial basis functions (RBF, Broomhead and Lowe, 1988), which is substantially cheaper (though still much more expensive than DE). This algorithm is referred to as RBFopt.

In addition, surrogate-based adaptive sampling is expected to be sensitive to dimensionality. For that reason this paper also considers a four dimensional version of the Branin-Hoo function, which is normally two dimensional.

## II. EGO, RBFopt and BFGS optimizers

**BFGS**

The SciPy (Virtanen et al. 2020) L-BFGS-B (Morales and Nocedal 2011) is a popular gradient based algorithm for local optimization (Venter 2010). The BFGS implementation uses finite differences to approximate the gradient of an objective function. In both of the following surrogate optimization algorithms, BFGS is run on the surrogate model. BFGS is used with EGO to find the maximum EI location, while BFGS is also used to find minimum locations of the RBF.

Global optimization algorithms, like EGO or RBFopt, are notorious for having poor local convergence. To address the poor local convergence, BFGS is sometimes run using the EGO optimum as the initial point. This strategy is referred to as EGO/BFGS in the paper. This is also done for RBFopt and referred to as RBFopt/BFGS.

**EGO**

The GPyOpt authors (2016) created GPyOpt which is the EGO implementation used in this work. The project is available online at https://github.com/SheffieldML/GPyOpt and has many rich features. The initial design is created using a random Latin hypercube sampling (LHS, Viana 2016). GPyOpt has many choices for surrogate models and acquisition functions, however we used the standard GP model with EI as the acquisition function which closely follows the original EGO description (Jones Et al. 1998). The algorithm is run for a set number of function evaluations.

**RBFopt**

The overhead of fitting a GP model and evaluating the maximum EI location with EGO can be substantial when the cost of the objective function is extremely cheap. A simple surrogate-based optimization algorithm was created for illustration purposes in this paper. The intent was to use a cheaper surrogate and acquisition function, such that more optimization iterations can be performed on the cheap objective function in the same amount of time (when compared to using GP with EI). The choice of surrogate is a multiquadric kernel ($\phi(r) = \sqrt{1 + r^2}$) RBF function, which should be cheaper to fit than a GP.

```
my_opt = sbopt.RbfOpt(my_fun,   # your objective function to minimize
                      bounds,   # bounds for your design variables
                      initial_design='latin',
                      initial_design_ndata=n_initial,   # No. of initial LHS points
                      n_local_optimze=5,
                      polish=polish,   # True or False, whether to run L-BFGS after
                      rbf_function='multiquadric ',
                      epsilon=None,
                      smooth=0.0,
                      metric='euclidean',
                      acquisition='rbf',
                      exploration='distance')
# run the optimizer
res = my_opt.minimize(verbose=0,
                      max_iter=max_iter,   # Depends on sample run
                      eps=0.002,
                      strategy='all_local',
                      n_same_best=20)
```
*Figure 1: Sample keywords used to drive the RBFopt algorithm.*

The optimization algorithm is available online at https://github.com/cjekel/sbopt. The specific keywords used in this study are presented in Figure 1. The algorithm can be described as the following:

1. Perform an initial design with LHS sampling and then fit a RBF.

2. Find the minimum of the RBF by performing multiple BFGS optimizations on the RBF. One optimization starts from best observed objective function value, while the remaining BFGS optimizations start from random points in the design space. The number of local BFGS runs is depicted by the *n_local_optimize* keyword. This begins an optimization iteration.

3. Using the *all_local* strategy, all of the minimum locations found on the RBF function are considered locations to evaluate the objective function. The objective function is evaluated at these locations, if and only if the Euclidean distance is at least a distance of *eps=0.002* from all previous design points. If none of the local minima satisfy this criterion, an optimization problem is solved to find a point which maximizes the minimum distance from all previous data points. This point would then be evaluated on the objective function, and this strategy is the *exploration='distance'* keyword. The newly evaluated responses are added to the database.

4. A new RBF is fit to the response using all of the evaluated design points. This completes an optimization iteration.

5. Steps 2-4 are repeated. Two stopping criteria are used to terminate the optimization process. The first considers whether the maximum number of iterations has been exceeded. The second criterion considers whether the best objective function value has improved in some fixed number of iterations (using the *n_same_best=20* keyword).

Preliminary runs have indicated that using about half the function evaluation budget for the initial sample is close to optimal. So in comparing the number of required function evaluations with the original and fortified Branin-Hoo functions we followed this principle. For EGO this was accomplished by specifying the number of iterations to be equal to the number of initial points. For the RBF optimizer it required playing with the number of initial points and number of local searches.

**Parallelization**

Multiple independent runs are performed to evaluate the average performance of the optimization algorithm on the test function. Several optimization runs were performed concurrently in parallel on a desktop computer to speed up the evaluation of the algorithm. This works well because each individual optimization run was effectively a single threaded task which did not benefit parallelization. The random seeds of the independent runs were controlled such that the entire evaluation could be replicated. The Joblib library available online at https://joblib.readthedocs.io was used to manage the parallel computing.

III. **The Branin-Hoo function**

Using the information from https://www.sfu.ca/~ssurjano/branin.html the Branin-Hoo function is defined as

$$b_r(\mathbf{x}) = b_r(x_1, x_2) = a\left(x_2 - bx_1^2 + cx_1 - r\right)^2 + s(1-t)\cos(x_1) + s.$$

We use the recommended values $a = 1, b = 5.1/(4\pi^2), c = 5/\pi, r = 6, s = 10, t = 1/(8\pi)$

The domain is here as usual, $x_1 \in [-5, 10], x_2 \in [0, 15]$, and the function is shown in Figure 2.

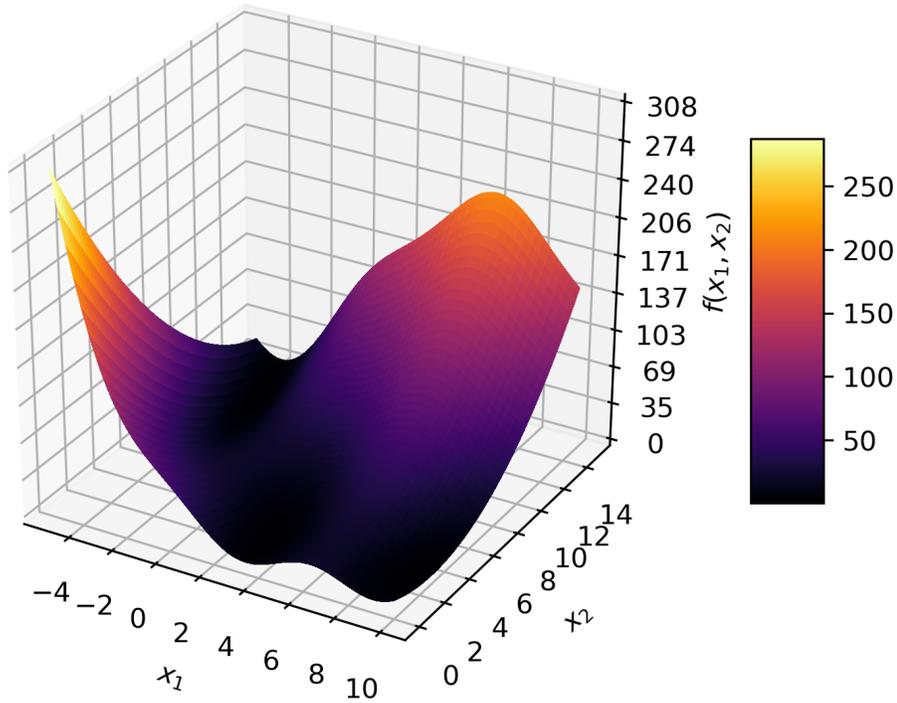

*Figure 2: The Branin-Hoo function.*

The function has three global optima where its value *is $b_r$ (x)=0.397887*. They are numbered here as

First: $\mathbf{x} = (-\pi, 12.275)$

Second: $\mathbf{x} = (\pi, 2.275)$

Third: $\mathbf{x} = (9.42478, 2.475)$

IV. **The 4D double Branin-Hoo function**

For the purpose of this study, we create a four-dimensional version of the Branin-Hoo function as

$$b_{r4}(\mathbf{x}) = b_{r4}(x_1, x_2, x_3, x_4) = b_r(x_1, x_2) + b_r(x_3, x_4).$$

This function has 9 equal-value optima. That is there are nine global optima, when $(x_1, x_2)$ take any of the values of the three optima of $b_r$ and $(x_3, x_4)$ take any value of the three optima. For example, the one combining the first and the second optima is $\mathbf{x} = (-\pi, 12.275, \pi, 2.275)$. At all nine optima $b_{r4} = 2 \times 0.397887 = 0.795774$.

V. **The bump**

It is desirable to add (for maximization) or subtract (for minimization) a bump to the original function that will not change the location of an optimum. Radial basis functions (RBF, Broomhead and Lowe, 1988) are selected because they depend only on the distance from the optimum. It is desirable that the bump will affect only one optimum, and for that the RBF[3] bump function is selected. The bump function is defined as

$$\varphi(r) = \begin{cases} \exp\left(-\dfrac{1}{1-(\varepsilon r)^2}\right) & \text{for } r < \dfrac{1}{\varepsilon} \\ 0 & \text{otherwise} \end{cases}$$

.

Here $r$ is the radial distance from the center of the bump, and its maximum value at $r=0$ is exp(-1)=0.3679. The width of the bump is determined by $\varepsilon$. Figure 3 shows a one-dimensional slice of the Branin-Hoo function with $10\varphi(r)$ bump subtracted at the location of its first global optimum.

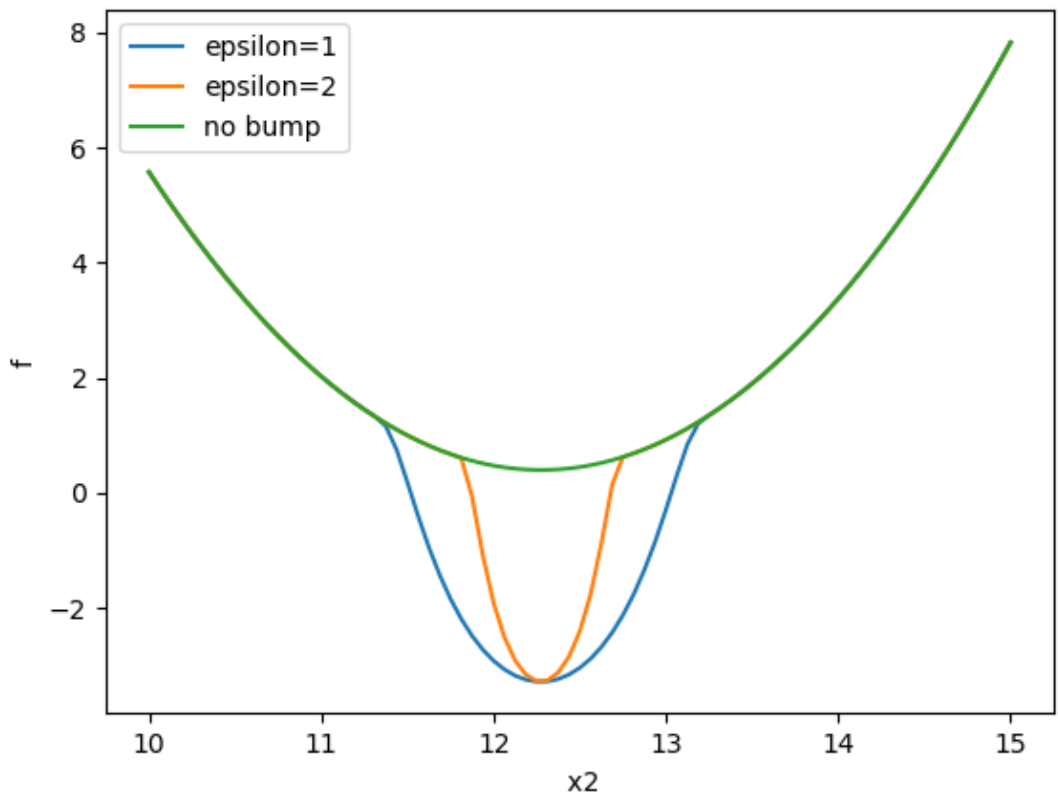

Figure 3: A 1D slice through the Branin-Hoo function for $x_1 = -\pi$ with the bump defined by Eq. 1.2 multiplied by 10 subtracted at the location of its first global optimum at $(-\pi, 12.275)$. The width parameters are $\varepsilon = 1$ and $\varepsilon = 2$.

---

[3] https://en.wikipedia.org/wiki/Radial_basis_function

VI. **Local and global optima of Branin-Hoo with a bump**

For the 2D Branin-Hoo function we apply a bump with an amplitude of 10 (subtracting 3.678794), as shown in Figure 3 to one of the optima. This leads to that optimum having a value of -3.280907, while the other two optima remain at 0.397887.

For the 4D Branin-Hoo function we apply a bump to $b_r(x_1, x_2)$ at optimum ($x_1 = -\pi, x_2 = 12.275$) and a bump to $b_r(x_3, x_4)$ at optimum ($x_3 = -\pi, x_4 = 12.275$). Both bumps use an amplitude of 5. Each bump subtracts 1.839397 from the function. So if a bump is applied to only one optimum, the function value is $0.795774 - 1.839397 = -1.043623$. If the bump is applied to both optima, the function value is $0.795774 - 2 \times 1.839397 = -2.88302$. These values of the local and global optima are shown in Table 1.

*Table 1: The optima of the 4D Branin-Hoo function with bumps applied to the first two variables at optimum b1 and to the last two variables at optimum b2.*

|  |  | $b_r(x_1, x_2)$ | | |
|---|---|---|---|---|
|  |  | $x_1 = -\pi, x_2 = 12.275$ | $x_1 = \pi, x_2 = 2.275$ | $x_1 = 9.42478, x_2 = 2.475$ |
| $b_r(x_3, x_4)$ | $x_3 = -\pi, x_4 = 12.275$ | -2.88302 (global) [b11] | -1.043623 [b21] | -1.043623 [b31] |
|  | $x_3 = \pi, x_4 = 2.275$ | -1.043623 [b12] | 0.795774 [b22] | 0.795774 [b32] |
|  | $x_3 = 9.42478, x_4 = 2.475$ | -1.043623 [b13] | 0.795774 [b33] | 0.795774 [b33] |

VII. **Performance measurement with single and double runs**

EGO and RBFopt are stochastic algorithms due to the random initial DOE, so that every time they are run one can get a different result. Therefore, when comparing the performance with and without a bump, we make multiple runs and compare the probability of failure deduced from the percentage of runs that fail to reach the optimum with some given tolerance. Here we need to choose the number of runs, which is designated as *n* needed for desired accuracy in the probability estimate.

The accuracy of the probability of failure estimate is measured here by the standard deviation of the number of runs that failed. It is easy to show that the standard deviation $\sigma_{fail}$ of the number of failures in *n* runs is related to the probability of failure *p* as

$$\sigma_{fail} = p\sqrt{pn(1-p)}.$$

To achieve accuracy of one percent, one sets $\sigma_{fail} = 0.01n$ and calculates the required number of runs for a given p.

$$(0.01n)^2 = p^3 n(1-p) \Rightarrow n = 10,000 p^3 (1-p).$$

While *p* is not known ahead of time, the maximum required *n* is attained for *p=0.75*, which gives *n=1055* runs. Since p=0.75 is the worst case, we rounded n to 1,000 for ease of translating the number of failures to percentages.

For some cases it was found that it is easy to get probability of failure on the order of 5%, but the number of function evaluations needed to get it below 1% is much higher. In that case, the probability of failing two independent runs $p_{double}$ is

$$p_{double} = p^2.$$

So as long as p is below 10%, $p_{double} \leq 1\%$.

## VIII. Performance with the 2D Branin-Hoo function

EGO and RBFopt were first run 1,000 times with the original 2D Branin-Ho function. Both were run with and without BFGS as a follow-up to improve the result.

Table 2 provides a summary of the best results obtained with 1,000 runs. It is seen that the combination of EGO and BFGS is extremely powerful, achieving perfect success with an average of 32 function evaluations. Without the BFGS follow-up, even 64 function evaluations were not enough for EGO. RBFopt required around 50 evaluations to achieve less than 1% probability of failure. Adding the BFGS follow-up would always result in a 0% probability of failure, as it appears BFGS is always successful at finding a local optimum from any point in the design space. It is also seen that both of these optimization algorithms have a preference for Optimum 2, which is in the middle of the design space compared to the other two optima.

*Table 2: Performance of algorithms for the original Branin-Hoo function without a bump. Results are based on 1000 replicate runs.*

| Algorithm | Initial number of points | Maximum Iterations | Percent failures | Average number of function evaluations | Percentage at or near each optimum B1, B2, B3 | | |
|---|---|---|---|---|---|---|---|
| EGO | 3 | 3 | 100% | 6 | 13 | 33 | 5 |
| EGO | 16 | 16 | 94% | 32 | 22 | 58 | 19 |
| EGO | 32 | 32 | 71% | 64 | 38 | 27 | 35 |
| EGO/BFGS | 3 | 3 | 0.0% | 32 | 29 | 46 | 25 |
| RBFopt/BFGS | 5 | 1 | 0.0% | 32.8 | 22 | 41 | 38 |
| RBFopt | 5 | 5 | 98.1% | 13.1 | 12 | 29 | 24 |
| RBFopt | 10 | 10 | 14.6% | 31.6 | 14 | 58 | 28 |
| RBFopt | 16 | 16 | 0.2% | 48.5 | 14 | 59 | 27 |
| RBFopt/BFGS | 16 | 16 | 0.0% | 63.5 | 14 | 59 | 27 |
| RBFopt | 25 | 1 | 98.5% | 27.7 | 19 | 37 | 22 |
| RBFopt | 25 | 25 | 0.0% | 67.4 | 25 | 51 | 24 |

### Performance with Bump at the first Optimum

The procedure described with no bump was repeated with the wider bump shown in Fig. 2, that is for an amplitude of *10/e*, with $\varepsilon = 1$.

The results when the bump was added to Optimum 1 are summarized in Table 3. With the larger number of function evaluations, we ran EGO for only 100 replicates, while RBFopt was continued with 1,000 replicates.

*Table 3: Performance of algorithms for the modified Branin-Hoo function with a bump amplitude of 10. Results for EGO are based on 100 replicate runs, but results for RBFopt are based on 1,000 replicates.*

| Algorithm | Initial number of points | Maximum Iterations | Percent failures | Average number of function evaluations | Percentage at or near each optimum B1, B2, B3 | | |
|---|---|---|---|---|---|---|---|
| EGO | 20 | 20 | 94% | 40 | 91 | 4 | 5 |
| EGO/BFGS | 20 | 20 | 9% | 61 | 91 | 4 | 5 |
| EGO | 25 | 25 | 89% | 50 | 100 | 0 | 0 |
| EGO/BFGS | 25 | 25 | 0.0% | 71 | 100 | 0 | 0 |
| EGO | 100 | 100 | 16% | 200 | 100 | 0 | 0 |
| RBFopt | 16 | 16 | 53% | 49.8 | 77 | 16 | 7 |
| RBFopt/BFGS | 16 | 16 | 23% | 67.1 | 77 | 16 | 7 |
| RBFopt | 25 | 25 | 12% | 70.1 | 98 | 1 | 0 |
| RBFopt/BFGS | 25 | 25 | 2% | 85.8 | 98 | 1 | 0 |
| RBFopt | 50 | 50 | 0.0% | 100.1 | 100 | 0 | 0 |
| RBFopt/BFGS | 50 | 50 | 0.0% | 114.9 | 100 | 0 | 0 |
| RBFopt | 80 | 20 | 0.8% | 113.8 | 100 | 0 | 0 |
| RBFopt | 90 | 10 | 13% | 112.4 | 99 | 0 | 0 |
| RBFopt | 100 | 5 | 42% | 113.1 | 92 | 1 | 0 |
| RBFopt | 100 | 10 | 10% | 122.0 | 99 | 0 | 0 |
| RBFopt | 100 | 20 | 0.1% | 132.7 | 100 | 0 | 0 |
| RBFopt | 200 | 5 | 16% | 212.5 | 99 | 0 | 0 |

It is seen that EGO and RBFopt required about doubled the number of required function evaluations to get similar performance as the original Branin-Hoo function. As before, it appears that the BFGS follow-up algorithm is successful at finding the global optima if started from the correct basin. Analyzing the first two rows of Table 3, it is seen that EGO brings 91% of the runs to the correct basin, but only 6% are close enough to the optimum. In contrast the 16-16 RBFopt which expended about 50 function evaluations only resulted in 77% of the runs end up in the correct basin, however a larger percentage of those runs were closer to the local optimum. The BFGS follow-up brings all the runs that ended up in the basin of Optimum 1 close enough to that optimum. Increasing the number of function evaluations slightly from 61 to 71, was sufficient with EGO/BFGS, while with EGO alone even 200 function evaluations were not enough. RBFopt in contrast was able to find the global optima in less than 101 function evaluations with 50 initial points and a 50 iteration limit.

IX. **Performance with no bump of the 4D Branin-Hoo function**

Without a bump there are nine equal-value optima. The best results are summarized in Table 4. We see that compared to Table 2, EGO/BFGS increased the number of function evaluations by almost a factor of 3, from 32 to 87, and similar with RBFopt/BFGS from 32 to 81.

EGO by itself has 100% failure even with 200 function evaluations, so it definitely benefits from using the BFGS follow-up. In contrast, RBFopt by itself can get down to low probability of failure.

It is curious that without the BFGS follow-up, RBFopt requires almost 400 function evaluations, and that progress towards 0% failure is not monotonic. This is due to the existence of nine basins. As points are added, RBFopt can model additional basins, but not necessarily accurately enough to get close enough to the optima in all 100 runs. Here one may consider doing two runs with 100 initial points and 100 iterations (row 17 in Table 4). With a failure rate of about 1% in a single run, the equation in VII leads us to expect a failure rate of about 0.01% in two runs, requiring 2x243=486 function evaluations.

*Table 4: Performance of algorithms for the 4D Branin-Hoo function without a bump. Results are based on 100 replicate runs.*

| Algorithm | Initial number of points | Maximum Iterations | Percent failures | Average number of function evaluations | Percentage at or near each optimum in the following order, b11, b12, b13, b21, b22, b23, b31, b32, b33 | | | | | | | | |
|---|---|---|---|---|---|---|---|---|---|---|---|---|---|
| EGO | 8 | 8 | 100% | 16 | 0 | 0 | 0 | 0 | 0 | 0 | 0 | 1 | 0 |
| EGO/BFGS | 8 | 8 | 0% | 84 | 7 | 10 | 12 | 13 | 16 | 12 | 2 | 12 | 16 |
| EGO | 50 | 50 | 100% | 100 | 2 | 4 | 2 | 6 | 15 | 7 | 3 | 5 | 3 |
| EGO | 100 | 100 | 100% | 200 | 11 | 25 | 4 | 17 | 32 | 4 | 0 | 2 | 0 |
| EGO | 250 | 2 | 100% | 252 | 4 | 0 | 2 | 1 | 1 | 3 | 1 | 2 | 0 |
| EGO | 250 | 10 | 100% | 260 | 4 | 0 | 4 | 1 | 2 | 2 | 5 | 1 | 0 |
| EGO | 250 | 25 | 100% | 275 | 8 | 6 | 12 | 7 | 6 | 3 | 4 | 3 | 5 |
| RBFopt | 10 | 1 | 100% | 11 | 0 | 0 | 1 | 0 | 0 | 1 | 0 | 0 | 1 |
| RBFopt/BFGS | 10 | 1 | 0% | 81 | 6 | 11 | 9 | 12 | 21 | 7 | 5 | 16 | 13 |
| RBFopt | 10 | 10 | 100% | 28 | 0 | 1 | 1 | 3 | 7 | 7 | 0 | 12 | 7 |
| RBFopt | 25 | 25 | 17% | 91 | 5 | 11 | 3 | 6 | 32 | 19 | 4 | 16 | 4 |
| RBFopt/BFGS | 25 | 25 | 0% | 133 | 5 | 11 | 3 | 6 | 32 | 19 | 4 | 16 | 4 |
| RBFopt | 50 | 2 | 100% | 56 | 0 | 1 | 1 | 2 | 3 | 2 | 2 | 4 | 0 |
| RBFopt | 50 | 50 | 0% | 175 | 3 | 14 | 4 | 12 | 41 | 11 | 2 | 1 | 3 |
| RBFopt/BFGS | 50 | 50 | 0% | 210 | 3 | 14 | 4 | 12 | 41 | 11 | 2 | 1 | 3 |
| RBFopt | 100 | 2 | 100% | 107 | 1 | 7 | 5 | 1 | 4 | 3 | 0 | 3 | 1 |
| RBFopt | 100 | 100 | 1% | 243 | 4 | 12 | 8 | 11 | 42 | 12 | 1 | 6 | 4 |
| RBFopt/BFGS | 100 | 100 | 0% | 279 | 4 | 12 | 8 | 11 | 42 | 12 | 1 | 6 | 4 |
| RBFopt | 250 | 2 | 100% | 258 | 5 | 8 | 2 | 6 | 10 | 4 | 3 | 2 | 1 |
| RBFopt | 250 | 10 | 97% | 289 | 10 | 16 | 6 | 11 | 27 | 8 | 4 | 8 | 4 |
| RBFopt | 250 | 25 | 5% | 342 | 7 | 15 | 4 | 12 | 38 | 9 | 4 | 7 | 4 |
| RBFopt | 300 | 2 | 100% | 308 | 3 | 8 | 3 | 4 | 16 | 5 | 5 | 2 | 2 |
| RBFopt | 300 | 30 | 0% | 407 | 10 | 23 | 3 | 13 | 32 | 9 | 3 | 6 | 1 |

Table 5: Performance of algorithms for the 4D Branin-Hoo function with a bump at (1,1). Results are based on 20 replicate runs for EGO, and 100 runs for RBFopt.

| Algorithm | Initial number of points | Maximum Iterations | Percent failures | Average number of function evaluations | Percentage at or near each optimum in the following order, b11, b12, b13, b21, b22, b23, b31, b32, b33 | | | | | | | | |
|---|---|---|---|---|---|---|---|---|---|---|---|---|---|
| EGO | 150 | 150 | 100% | 300 | 75 | 5 | 5 | 15 | 0 | 0 | 0 | 0 | 0 |
| EGO/BFGS | 150 | 150 | 25% | 353 | 75 | 10 | 5 | 10 | 0 | 0 | 0 | 0 | 0 |
| EGO | 250 | 250 | 100% | 500 | 100 | 0 | 0 | 0 | 0 | 0 | 0 | 0 | 0 |
| EGO/BFGS | 250 | 250 | 0% | 552 | 100 | 0 | 0 | 0 | 0 | 0 | 0 | 0 | 0 |
| RBFopt | 150 | 150 | 41% | 347 | 61 | 10 | 2 | 23 | 0 | 0 | 4 | 0 | 0 |
| RBFopt/BFGS | 150 | 150 | 39% | 383 | 61 | 10 | 2 | 23 | 0 | 0 | 4 | 0 | 0 |
| RBFopt | 250 | 250 | 36% | 443 | 66 | 14 | 4 | 14 | 0 | 0 | 2 | 0 | 0 |
| RBFopt/BFGS | 250 | 250 | 34% | 477 | 66 | 14 | 4 | 14 | 0 | 0 | 2 | 0 | 0 |
| RBFopt | 300 | 10 | 99% | 340 | 25 | 2 | 8 | 18 | 8 | 4 | 9 | 2 | 2 |
| RBFopt | 500 | 8 | 98% | 532 | 32 | 25 | 4 | 17 | 4 | 1 | 6 | 0 | 0 |
| RBFopt | 500 | 500 | 17% | 678 | 88 | 7 | 0 | 5 | 0 | 0 | 0 | 0 | 0 |
| RBFopt/BFGS | 500 | 500 | 12% | 711 | 88 | 7 | 0 | 5 | 0 | 0 | 0 | 0 | 0 |
| RBFopt | 700 | 9 | 94% | 736 | 41 | 25 | 4 | 21 | 2 | 0 | 3 | 0 | 0 |
| RBFopt/BFGS | 700 | 9 | 55% | 789 | 45 | 25 | 4 | 21 | 2 | 0 | 3 | 0 | 0 |
| RBFopt | 1000 | 14 | 71% | 1054 | 72 | 14 | 3 | 8 | 0 | 0 | 3 | 0 | 0 |
| RBFopt/BFGS | 1000 | 14 | 28% | 1099 | 72 | 14 | 3 | 8 | 0 | 0 | 3 | 0 | 0 |
| RBFopt | 1000 | 20 | 51% | 1075 | 84 | 8 | 0 | 7 | 0 | 0 | 1 | 0 | 0 |
| RBFopt/BFGS | 1000 | 20 | 16% | 1115 | 84 | 8 | 0 | 7 | 0 | 0 | 1 | 0 | 0 |
| RBFopt | 1000 | 1000 | 9% | 1145 | 97 | 1 | 0 | 2 | 0 | 0 | 0 | 0 | 0 |
| RBFopt/BFGS | 1000 | 1000 | 3% | 1176 | 97 | 1 | 0 | 2 | 0 | 0 | 0 | 0 | 0 |

X. **Performance of bump at (1,1) for 4D Branin-Hoo function**

Table 5 shows selected cases with the bump at (1,1). With the large numbers of points required, we could afford only 20 replicate EGO runs, and 100 runs for RBFopt. We again see that EGO absolutely needs the BFGS follow-up, while RBFopt can get relatively low probabilities of failure without the polish. However, it should be noted that EGO is better at identifying the correct basin of the global optimum. EGO with 250-250 resulted with 20/20 runs finding the correct basin, while RBFopt with 500-500 resulted in 88/100 runs finding the correct basin.

With the BFGS follow-up, the performance of EGO is impressive in terms of function evaluations compared to the less expensive RBFopt/BFGS. It needs only about 550 function evaluations compared to about 1200 for RBFopt/BFGS. All of the RBFopt optimization runs with 1000-1000 did not reach the maximum number of iterations as they were stopped early after observing 20 iterations without improvement.

Compared to the unfortified version of the four-dimensional Branin-Ho (compare to Table 4), EGO and RBFopt required more than 6 times the number of function evaluations to find the correct basin. That is the fortified 4D function is much more of a challenge to both algorithms compared to the 2D function. This reflects the effect of dimensionality and the larger number of basins.

## XI. Concluding Remarks

Global optimization algorithms are often tested on easy functions with multiple identical-valued optima. This paper suggests that these functions could be fortified to become much harder by adding or subtracting a bump to one of the optima. This was illustrated in a previous paper by the authors for the Branin-Hoo function, which has three identical-valued global optima. In that previous work the performance of the Python SciPy popular differential evolution (DE) optimizer with a follow up (called polish) by the BFGS local gradient based algorithm. In this paper we examined the performance of two surrogate-based optimization algorithms that proceed by adaptive sampling. The first algorithm is the well-known EGO algorithm, which is based on a Gaussian Process (GP) surrogate. The second algorithm is a much cheaper algorithm based on radial basis function denoted RBFopt that allowed us to run more replicate runs needed to achieve more accurate estimates of the probability of failing to reach the global optimum. These algorithms were applied to the original and fortified Branin-Hoo function and also to a four dimensional version of the Branin-Hoo function that has nine identical-valued global optima with and without fortification. The effect of fortification on the number of function evaluations led us to the following conclusions:

1. The fortification was associated with a factor of 2 or 3 increase in the number of function evaluations for the 2D Branin-Hoo function and 7-14 increase for the 4D function.
2. The fortification can change the relative performance of two competing algorithms. For the 4D case, without fortification RBFopt was able to find the local optima with fewer function evaluations, however with fortification EGO was significantly better at identifying the correct basin containing the global optima.
3. Surrogate-based optimization is sensitive to dimensionality, so extensions of the fortified Branin-Hoo function to higher dimensions may be useful for evaluating surrogate-based optimization algorithms.
4. Jekel and Haftka (2020) showed that DE required 16 times more function evaluations with the fortified 2D Branin-hoo function. However, the surrogate-based optimization algorithms presented here required about only double the number of function evaluations, demonstrating the frugality with this class of optimization algorithm.

**References**

Broomhead, D. S., & Lowe, D. (1988). Radial basis functions, multi-variable functional interpolation and adaptive networks, (No. RSRE-MEMO-4148). Royal Signals and Radar Establishment Malvern (United Kingdom).

Forrester, A. I., & Keane, A. J. (2009). Recent advances in surrogate-based optimization. Progress in aerospace sciences, 45(1-3), 50-79.

Jekel, C. F., Haftka, R. T. (2019). Fortified Test Functions for Global Optimization and the Power of Multiple Runs. arXiv preprint arXiv:1912.10575. https://arxiv.org/abs/1912.10575

Jekel, C. F., Haftka, R. (2020). Weaponizing Favorite Test Functions for Testing Global Optimization Algorithms: An illustration with the Branin-Hoo Function. In AIAA AVIATION 2020 FORUM (p. 3132).

Jones, D. R., Schonlau, M., & Welch, W. J. (1998). Efficient global optimization of expensive black-box functions. Journal of Global optimization, 13(4), 455-492.

Morales, J.L. Nocedal, J. (2011) L-BFGS-B: Remark on Algorithm 778: L-BFGS-B, FORTRAN routines for large scale bound constrained optimization, ACM Transactions on Mathematical Software, 38, 1.

Storn, R., & Price, K. (1997). Differential Evolution -- A Simple and Efficient Heuristic for global Optimization over Continuous Spaces. Journal of Global Optimization, 11(4), 341–359. https://doi.org/10.1023/A:1008202821328

The GPyOpt authors. (2016). GPyOpt: a Bayesian optimization framework in Python.

Venter, G. (2010). Review of optimization techniques.

Viana, F. A. (2016). A tutorial on Latin hypercube design of experiments. Quality and reliability engineering international, 32(5), 1975-1985.

Virtanen, P., Gommers, R., Oliphant, T.E. *et al.* (2020) SciPy 1.0: fundamental algorithms for scientific computing in Python. *Nat Methods*. https://doi.org/10.1038/s41592-019-0686-2